\documentclass[12pt]{article}      
\usepackage{amssymb}
\usepackage{amscd}

\input xy       
\xyoption{all}      

\newtheorem{theorem}{Theorem}[section]
\newtheorem{prop}[theorem]{Proposition}

\newtheorem{lemma}[theorem]{Lemma}
\newtheorem{coro}[theorem]{Corollary}
\newtheorem{prop-def}{Proposition-Definition}[section]

\newcommand{\nc}{\newcommand}

\nc{\bin}[2]{ (_{\stackrel{\scs{#1}}{\scs{#2}}})}  
\nc{\binc}[2]{(\!\! \begin{array}{c} \scs{#1}\\
    \scs{#2} \end{array}\!\!)}  
\nc{\bincc}[2]{  ( {\scs{#1} \atop
    \vspace{-1cm}\scs{#2}} )}  
\nc{\bs}{\bar{S}}
\nc{\la}{\longrightarrow}
\nc{\rar}{\rightarrow}
\nc{\dar}{\downarrow}
\nc{\dap}[1]{\downarrow \rlap{$\scriptstyle{#1}$}}
\nc{\defeq}{\stackrel{\rm def}{=}}
\nc{\dis}[1]{\displaystyle{#1}}
\nc{\dotcup}{\ \displaystyle{\bigcup^\bullet}\ }
\nc{\hcm}{\ \hat{,}\ }
\nc{\hts}{\hat{\otimes}}
\nc{\hcirc}{\hat{\circ}}
\nc{\lleft}{[}
\nc{\lright}{]}
\nc{\curlyl}{\left \{ \begin{array}{c} {} \\ {} \end{array}
    \right .  \!\!\!\!\!\!\!}
\nc{\curlyr}{ \!\!\!\!\!\!\!
    \left . \begin{array}{c} {} \\ {} \end{array}
    \right \} }
\nc{\longmid}{\left | \begin{array}{c} {} \\ {} \end{array}
    \right . \!\!\!\!\!\!\!}
\nc{\ora}[1]{\stackrel{#1}{\rar}}
\nc{\ola}[1]{\stackrel{#1}{\la}}
\nc{\scs}[1]{\scriptstyle{#1}}
\nc{\mrm}[1]{{\rm #1}}
\nc{\dirlim}{\displaystyle{\lim_{\longrightarrow}}\,}
\nc{\invlim}{\displaystyle{\lim_{\longleftarrow}}\,}
\nc{\dislim}[1]{\displaystyle{\lim_{#1}}}
\nc{\colim}{\mrm{colim}}
\nc{\mvp}{\vspace{0.3cm}}
\nc{\tk}{^{(k)}}
\nc{\tp}{^\prime}
\nc{\ttp}{^{\prime\prime}}
\nc{\svp}{\vspace{2cm}}
\nc{\vp}{\vspace{8cm}}
\nc{\proofend}{$\blacksquare$ \vspace{0.3cm}}
\nc{\modg}[1]{\!<\!\!{#1}\!\!>}
\nc{\intg}[1]{F_C(#1)}
\nc{\lmodg}{\!<\!\!}
\nc{\rmodg}{\!\!>\!}
\nc{\cpi}{\widehat{\Pi}}
\nc{\sha}{{\mbox{\cyr X}}}  
\nc{\shpr}{\diamond}    
\nc{\labs}{\mid\!}
\nc{\rabs}{\!\mid}

\nc{\ann}{\mrm{ann}}
\nc{\Aut}{\mrm{Aut}}
\nc{\can}{\mrm{can}}
\nc{\Cont}{\mrm{Cont}}
\nc{\rchar}{\mrm{char}}
\nc{\cok}{\mrm{coker}}
\nc{\dtf}{{R-{\rm tf}}}
\nc{\dtor}{{R-{\rm tor}}}

\nc{\Div}{{\mrm Div}}
\nc{\End}{\mrm{End}}
\nc{\Ext}{\mrm{Ext}}
\nc{\Fil}{\mrm{Fil}}
\nc{\Fr}{\mrm{Fr}}
\nc{\Frob}{\mrm{Frob}}
\nc{\Gal}{\mrm{Gal}}
\nc{\GL}{\mrm{GL}}
\nc{\Hom}{\mrm{Hom}}
\nc{\hsr}{\mrm{H}}
\nc{\hpol}{\mrm{HP}}
\nc{\id}{\mrm{id}}
\nc{\im}{\mrm{im}}
\nc{\incl}{\mrm{incl}}
\nc{\length}{\mrm{length}}
\nc{\mchar}{\rm char}
\nc{\mpart}{\mrm{part}}
\nc{\ql}{{\QQ_\ell}}
\nc{\qp}{{\QQ_p}}
\nc{\rank}{\mrm{rank}}
\nc{\rcot}{\mrm{cot}}
\nc{\rdef}{\mrm{def}}
\nc{\rdiv}{{\rm div}}
\nc{\rtf}{{\rm tf}}
\nc{\rtor}{{\rm tor}}
\nc{\res}{\mrm{res}}
\nc{\SL}{\mrm{SL}}
\nc{\Spec}{\mrm{Spec}}
\nc{\tor}{\mrm{tor}}
\nc{\Tr}{\mrm{Tr}}
\nc{\tr}{\mrm{tr}}

\nc{\ab}{\mathbf{Ab}}
\nc{\Alg}{\mathbf{Alg}}
\nc{\Bax}{\mathbf{Bax}}
\nc{\bfk}{{\bf k}}
\nc{\bfone}{{\bf 1}}
\nc{\detail}{\marginpar{\bf More detail}
    \noindent{\bf Need more detail!}
    \svp}
\nc{\Diff}{\mathbf{Diff}}
\nc{\gap}{\marginpar{\bf Incomplete}\noindent{\bf Incomplete!!}
    \svp}
\nc{\FMod}{\mathbf{FMod}}
\nc{\Int}{\mathbf{Int}}
\nc{\Mon}{\mathbf{Mon}}
\nc{\margin}[1]{\marginpar{}}
\nc{\proof}{\noindent{\bf Proof: }}
\nc{\remarks}{\noindent{\bf Remarks: }}
\nc{\Rep}{\mathbf{Rep}}
\nc{\Rings}{\mathbf{Rings}}
\nc{\Sets}{\mathbf{Sets}}
\nc{\bill}[1]{\marginpar{\bf To Bill}\noindent{\bf To Bill:}
    {\tt #1}\\ }
\nc{\li}[1]{\marginpar{\bf To Li}\noindent{\bf To Li:}
    {\tt #1}\\ }

\nc{\BA}{{\Bbb A}}
\nc{\CC}{{\Bbb C}}
\nc{\DD}{{\Bbb D}}
\nc{\EE}{{\Bbb E}}
\nc{\FF}{{\Bbb F}}
\nc{\GG}{{\Bbb G}}
\nc{\HH}{{\Bbb H}}
\nc{\LL}{{\Bbb L}}
\nc{\NN}{{\Bbb N}}
\nc{\QQ}{{\Bbb Q}}
\nc{\RR}{{\Bbb R}}
\nc{\TT}{{\Bbb T}}
\nc{\VV}{{\Bbb V}}
\nc{\ZZ}{{\Bbb Z}}

\nc{\cala}{{\cal A}}
\nc{\calc}{{\cal C}}
\nc{\cald}{\mathcal{D}}
\nc{\cale}{{\cal E}}
\nc{\calf}{{\cal F}}
\nc{\calg}{{\cal G}}
\nc{\calh}{{\cal H}}
\nc{\cali}{{\cal I}}
\nc{\call}{{\cal L}}
\nc{\calm}{{\cal M}}
\nc{\caln}{{\cal N}}
\nc{\calo}{{\cal O}}
\nc{\calp}{{\cal P}}
\nc{\calr}{{\cal R}}
\nc{\cals}{{\cal S}}
\nc{\calt}{{\cal T}}
\nc{\calw}{{\cal W}}
\nc{\calx}{{\cal X}}
\nc{\CA}{\mathcal{A}}

\nc{\fraka}{{\frak a}}
\nc{\frakB}{{\frak B}}
\nc{\frakm}{{\frak m}}
\nc{\frakp}{{\frak p}}
\nc{\frakS}{{\frak S}}
\nc{\frakA}{{\frak A}}
\nc{\frakx}{{\frak x}}

\font\cyr=wncyr10

\title{
Properties of Free Baxter Algebras
\thanks{The author is supported in part by NSF grant
   \#DMS 97-96122.
  MSC Numbers: Primary 16A06, 47B99.
  Secondary 13A99,16W99.} }
\author{
Li Guo
\\
Department of Mathematics and Computer Science\\ Rutgers
University\\ Newark, NJ 07102, USA \\ (liguo@newark.rutgers.edu) }
\date{} 

\begin{document}
\maketitle

\vspace{.5cm}

\begin{center}
{\Large\bf Abstract}
\end{center}
The study of free Baxter algebras was started by Rota and Cartier
thirty years ago. We continue this study by applying two recent
constructions of free Baxter algebras. We investigate the basic
structure of a free Baxter algebra, and characterize in detail
when a free Baxter algebra is a domain or a reduced algebra. We
also describe the nilpotent radical of a free Baxter algebra when
it is not reduced.

\vspace{0.5cm}

\setcounter{section}{0}

\section{Introduction}

The study of Baxter operators originated in the work of
Baxter~\cite{Ba} on fluctuation theory, and the algebraic study of
Baxter operators was started by Rota~\cite{Ro1}. Let $C$ be a
commutative ring and let $\lambda$ be a fixed element in $C$. A
Baxter algebra of weight $\lambda$ is a commutative $C$-algebra
$R$ together with a $C$-linear operator $P$ on $R$ such that for
any $x,\ y\in R$,
\[ P(x)P(y) =P(xP(y))+P(yP(x))+\lambda P(xy).\]
Baxter algebras have important applications in
combinatorics~\cite{Ro2,Ro3} and are closely related to several
areas in algebra and geometry, such as differential
algebras~\cite{Ko}, difference algebras~\cite{R.Co} and iterate
integrals in geometry~\cite{Ch}.

As in any algebraic system, free Baxter algebras play a central
role in the study of Baxter algebras. Even though the existence of
free Baxter algebras follows from the general theory of universal
algebras, in order to get a good understanding of free Baxter
algebras, it is desirable to find concrete constructions of a free
Baxter algebra. Two constructions were given in~\cite{G-K1,G-K2},
called shuffle Baxter algebras and standard Baxter algebras
respectively (see Section~\ref{sec:back} for details). The
construction of shuffle Baxter algebras is motivated by the
shuffle product of iterated integrals~\cite{Re} and an earlier
construction of Cartier~\cite{Ca}. The construction of standard
Baxter algebras is motivated by a construction of Rota~\cite{Ro1}.

In this paper, we apply these two constructions of free Baxter
algebras to obtain further information about free Baxter algebras.
After a brief discussion of basic properties of free Baxter
algebras, we will focus on the investigation of zero divisors and
nilpotent elements in a free Baxter algebra. This question has
been considered by Cartier~\cite{Ca} and Rota~\cite{Ro1,Ro2} for
Baxter algebras of weight one without an identity. In their case,
the free Baxter algebras have very good properties. In fact the
algebras are often isomorphic to either polynomial algebras or
power series algebras. The explicit descriptions of free Baxter
algebras obtained in \cite{G-K1,G-K2} enable us to consider this
question for a more general class of Baxter algebras. It is
interesting to observe that even if a free Baxter algebra is
constructed from an integral domain or a reduced algebra, the free
Baxter algebra is not necessarily a domain or a reduced algebra.
We show that the obstruction
depends on several factors, including the characteristic of the
base algebra, the weight of the Baxter algebra, whether or not the
Baxter algebra has an identity and whether or not the Baxter
algebra is complete. We provide necessary and sufficient
conditions for a free Baxter algebra to be a domain or to be
reduced (Theorem~\ref{thm:domain} and \ref{thm:domain1}), and
describe the nilpotent radical when a free Baxter algebra is not
reduced (Theorem~\ref{thm:red}).

We first give a brief summary of
the concept of Baxter algebras and the two
constructions of free Baxter algebras in section~\ref{sec:back}.
In section~\ref{sec:basic} we study basic properties of
free Baxter algebras, such as subalgebras, quotient algebras
and limits.
In section~\ref{sec:domain} we study in detail when a free
Baxter algebra is a domain or a reduced algebra.
We also consider free complete Baxter algebras.

\section{Free Baxter algebras}
\margin{sec:back} \label{sec:back} For later application, we will
describe the constructions of free Baxter
algebras~\cite{G-K1,G-K2}. We will also prove some preliminary
results.

We write $\NN$ for the additive monoid of
natural numbers and
$\NN_+$ for the positive integers.
Any ring $C$ is commutative with
identity element $\bfone_C$, and any ring homomorphism
preserves the identity elements.
For any $C$-modules $M$ and $N$,
the tensor product $M\otimes N$ is taken over $C$
unless otherwise indicated.
For a $C$-module $M$ and
$n\in \NN_+$, denote the tensor power
\[ M^{\otimes n}=\underbrace{M\otimes\ldots\otimes M}
    _{n\ {\rm factors}}.\]

\subsection{Baxter algebras}
For a given ring $C$, let $\Alg_C$ denote the category of
commutative $C$-algebras with an identity. For a given $\lambda\in
C$ and $R\in \Alg_C$,
\begin{itemize}
\item
a {\bf Baxter operator of weight $\lambda$ on $R$ over $C$}
is a $C$-module
endomorphism $P$ of $R$ satisfying
\begin{equation}
 P(x)P(y)=P(xP(y))+P(yP(x))+\lambda P(xy),\ x,\ y\in R;
\label{eq:bax1}
\end{equation}
\item
a {\bf Baxter C-algebra of weight $\lambda$}
is a pair $(R,P)$ where $R$ is a $C$-algebra
and $P$ is a Baxter operator of weight $\lambda$
on $R$ over $C$.
\item
a $C$-algebra homomorphism
$f:R\to S$ between two Baxter $C$-algebras
$(R,P)$ and $(S,Q)$ of weight $\lambda$
is called a {\bf homomorphism of Baxter $C$-algebras}
if $ f(P(x))=Q(f(x))$ for all $x \in R$.
\end{itemize}
Denote $\Bax_{C,\lambda}$ for the category of Baxter $C$-algebras
of weight $\lambda$.
If the meaning of $\lambda$ is clear, we will
suppress $\lambda$ from the notation.

A {\bf Baxter ideal} of $(R,P)$ is an ideal $I$ of $R$ such that
$P(I)\subseteq I$. The concepts of Baxter subalgebras, quotient
Baxter algebras can be similarly defined. It follows from the
general theory of universal algebras that limits and colimits
exist in $\Bax_C$~\cite{P.Co},\cite[p84]{Ja}, \cite[p.210]{Ma}. In
particular, inverse limits and direct limits exist in $\Bax_C$.

\subsection{Shuffle Baxter algebras}
\margin{sec:shuf}
\label{sec:shuf}

For $m,n\in \NN_+$,
define the set of {\bf $(m,n)$-shuffles} by
\[ S(m,n)=
 \left \{ \sigma\in S_{m+n}
    \begin{array}{ll} {} \\ {} \end{array} \right .
\left |
\begin{array}{l}
\sigma^{-1}(1)<\sigma^{-1}(2)<\ldots<\sigma^{-1}(m),\\
\sigma^{-1}(m+1)<\sigma^{-1}(m+2)<\ldots<\sigma^{-1}(m+n)
\end{array}
\right \}.\]
Given an $(m,n)$-shuffle $\sigma\in S(m,n)$,
a pair of indices $(k, k+1)$,\ $1\leq k< m+n$ is
called an {\bf admissible pair} for $\sigma$
if $\sigma(k)\leq m<\sigma(k+1)$.
Denote $\calt^\sigma$ for the set of admissible pairs for $\sigma$.
For a subset $T$ of $\calt^\sigma$, call the pair
$(\sigma,T)$ a {\bf mixable $(m,n)$-shuffle}.
Let $\mid T\mid$ be the cardinality of $T$.
$(\sigma,T)$ is identified with $\sigma$
if $T$ is the empty set.
Denote
\[ \bs (m,n)=\{ (\sigma,T)\mid \sigma\in S(m,n),\
    T\subset \calt^\sigma\} \]
for the set of {\bf $(m,n)$-mixable shuffles}.

For $A\in \Alg_C$,
$x=x_1\otimes\ldots\otimes x_m\in A^{\otimes m}$,
$y=y_1\otimes \ldots\otimes y_n\in A^{\otimes n}$
and $(\sigma,T)\in \bar{S}(m,n)$,
the element
\[  \sigma (x\otimes y) =u_{\sigma(1)}\otimes u_{\sigma(2)} \otimes
    \ldots \otimes u_{\sigma(m+n)}\in A^{\otimes (m+n)},\]
where
\[ u_k=\left \{ \begin{array}{ll}
    x_k,& 1\leq k\leq m,\\
    y_{k-m}, & m+1\leq k\leq m+n, \end{array}
    \right . \]
is called a {\bf shuffle} of $x$ and $y$;
the element
\[ \sigma(x\otimes y; T)= u_{\sigma(1)}\hts u_{\sigma(2)} \hts
    \ldots \hts u_{\sigma(m+n)} \in A^{\otimes (m+n-\mid T\mid)},\]
where for each pair $(k,k+1)$, $1\leq k< m+n$,
\[ u_{\sigma(k)}\hts u_{\sigma(k+1)} =\left \{\begin{array}{ll}
    u_{\sigma(k)} u_{\sigma(k+1)},  &
     (k,k+1)\in T\\
    u_{\sigma(k)}\otimes u_{\sigma(k+1)}, &
    (k,k+1) \not\in T,
    \end{array} \right . \]
is called a {\bf mixable shuffle} of $x$ and $y$.

Fix a $\lambda\in C$.
Let

\[ \sha_C(A)=\sha_{C,\lambda}(A)= \bigoplus_{k\in\NN}
    A^{\otimes (k+1)}
= A\oplus A^{\otimes 2}\oplus \ldots \]
be the Baxter $C$-algebra of weight $\lambda$~\cite{G-K1} in which

\begin{itemize}
\item
the $C$-module structure
is the natural one,
\item
the multiplication is the {\bf mixed shuffle product},
defined by

\margin{eq:shuf}
\begin{equation}
x\shpr y\
=\sum_{(\sigma,T)\in \bs (m,n)} \lambda^{\mid T\mid }
    x_0y_0\otimes
\sigma(x^+\otimes y^+;T)
 \in \bigoplus_{k\leq m+n+1} A^{\otimes k}
\label{eq:shuf}
\end{equation}
for $x=x_0\otimes x_1\otimes\ldots \otimes x_m\in
A^{\otimes (m+1)}$ and
$y=y_0\otimes y_1\otimes\ldots\otimes y_n\in A^{\otimes (m+1)}$,
where $x^+=x_1\otimes \ldots \otimes x_m$ and
$y^+=y_1\otimes \ldots \otimes y_n$,
\item
the weight $\lambda$ Baxter operator $P_A$ on
$\sha_C(A)$ is obtained by assigning
\[ P_A( x_0\otimes x_1\otimes \ldots \otimes x_n)
=\bfone_A\otimes x_0\otimes x_1\otimes \ldots\otimes x_n, \]
for all
$x_0\otimes x_1\otimes \ldots\otimes x_n\in A^{\otimes (n+1)}$.
\end{itemize}
$(\sha_C(A),P_A)$ is called the {\bf shuffle Baxter $C$-algebra on
$A$ of weight $\lambda$}. When there is no danger of confusion, we
often suppress $\shpr$ in the mixed shuffle product.  To
distinguish the $C$-submodule $A^{\otimes k}$ of $\sha_C(A)$ from
the tensor power $C$-algebra $A^{\otimes k}$, we sometimes denote
$\sha_C^{k-1}(A)$ for $A^{\otimes k}\subseteq \sha_C(A)$.

For a given set $X$, we also let $(\sha_C(X),P_X)$ denote the
shuffle Baxter $C$-algebra $(\sha_C(C[X]),P_{C[X]})$, called the
{\bf shuffle Baxter $C$-algebra on $X$ (of weight $\lambda$).} Let
$j_A:A\rar \sha_C(A)$ (resp. $j_X:X\to \sha_C(X)$) be the
canonical inclusion map.

\begin{theorem}
\margin{thm:shua} \label{thm:shua} \cite{Ca,G-K1}
$(\sha_C(A),P_A)$, together with the natural embedding $j_A$, is a
free Baxter $C$-algebra on $A$ of weight $\lambda$. In other
words, for any Baxter $C$-algebra $(R,P)$ and any $C$-algebra
homomorphism $\varphi:A\rar R$, there exists a unique Baxter
$C$-algebra homomorphism $\tilde{\varphi}:(\sha_C(A),P_A)\rar
(R,P)$ such that the  diagram
\[\xymatrix{
A \ar[rr]^(0.4){j_A} \ar[drr]_{\varphi}
    && \sha_C(A) \ar[d]^{\tilde{\varphi}} \\
&& R } \]
commutes.
Further, any $f:A\to B$ in $\Alg_C$ extends uniquely
to
\[ \sha_C(f): \sha_C(A)\to \sha_C(B) \]
in $\Bax_C$.
More precisely,
$ \sha_C(f)= \oplus_{n\in \NN} f^{\otimes (n+1)}$
with
$f^{\otimes (n+1)}:A^{\otimes (n+1)}\to B^{(n+1)}$
being the $(n+1)$-th tensor power of the $C$-module
homomorphism $f$.
Similarly,
$(\sha_C(X),P_X)$, together with the natural embedding
$j_X$,  is a free Baxter $C$-algebra on $X$
of weight $\lambda$.
\end{theorem}

Taking $A=C$, we get
\[ \sha_C(C)=\bigoplus_{n=0}^\infty C^{\otimes (n+1)}
= C \bfone^{\otimes (n+1)}. \]
where
$\bfone^{\otimes (n+1)}
= \underbrace{\bfone_C \otimes \ldots \otimes \bfone_C}
_{(n+1)-{\rm factors}}$. In this case the mixable shuffle
product formula~(\ref{eq:shuf}) gives

\begin{prop}
\margin{prop:unit}
\label{prop:unit}
For any $m,n\in \NN$,
\[ \bfone^{\otimes (m+1)} \shpr \bfone^{\otimes (n+1)} =
\sum_{k=0}^m \binc{m+n-k}{n}\binc{n}{k} \lambda^k
\bfone^{\otimes (m+n+1-k)}.\]
\end{prop}

\subsection{Complete shuffle Baxter algebras}
We now consider the completion of $\sha_C(A)$.
Recall that we denote $\sha^k_C(A)$ for the $C$-submodule
$A^{\otimes (k+1)}$ of $\sha_C(A)$.

Given $k\in \NN_+$,
$\Fil^k \sha_C(A)\defeq \bigoplus_{n\geq k} \sha^n_C(A),$
is a Baxter ideal of $\sha_C(A)$.
Denote $\widehat{\sha}_C(A)= \invlim \sha_C(A)/\Fil^k\sha_C(A)$,
called the {\bf complete shuffle Baxter algebra on $A$},
with the Baxter operator denoted by $\hat{P}$.
It naturally contains $\sha_C(A)$ as a Baxter subalgebra
and is a free object in the category of Baxter algebras
that are complete with respect to a canonical filtration
defined by the Baxter operator~\cite{G-K2}.
On other hand, consider the infinite product of $C$-modules
$\prod_{k\in \NN} \sha_C^k (A)$.
It contains $\sha_C(A)$ as a dense subset with respect to the
topology defined by the filtration $\Fil^k \sha_C(A)$.
All operations of the Baxter $C$-algebra $\sha_C(A)$ are continuous
with respect to this topology, hence extend uniquely to
operations on $\prod_{k\in\NN} \sha_C^k(A)$,
making $\prod_{k\in\NN}\sha_C^k(A)$ a Baxter algebra
of weight $\lambda$, with the Baxter operator denoted by $\bar{P}$.

\begin{theorem}
\margin{thm:series}
\label{thm:series}
\cite{G-K2}
\begin{enumerate}
\item
The map
\[ \psi_{A}:
\widehat{\sha}_C(A) \to \prod_{k\in \NN} \sha_C^k (A),\
 ((x^{(n)}_k)_k + \Fil^n \sha_C(A))_n \mapsto
(x^{(k)}_k)_k \]
is an isomorphism of Baxter algebras
extending the identity map on $\sha_C(A)$.
\item
Given a morphism $f:A\to B$ in $\Alg_C$, we have the following
commutative diagram
\[ \begin{array}{ccc}
\widehat{\sha}_C(A) & \ola{\psi_A} &
    \prod_{k\in\NN} \sha_C^k (A) \\
\dap{\widehat{\sha}_C(f)} &&    \dap{\prod_k f_k} \\
\widehat{\sha}_C(B) & \ola{\psi_B} &
    \prod_{k\in\NN} \sha_C^k (B)
\end{array} \]
where $\widehat{\sha}_C(f)$ is induced from
$\sha_C(f)$ in Theorem~\ref{thm:shua} by taking completion,
and $f_k:\sha_C^k(A) \to \sha_C^k(B)$ is the
tensor power morphism of $C$-modules
$f^{\otimes (k+1)}: A^{\otimes (k+1)} \to B^{\otimes (k+1)}$
induced from $f$.
\end{enumerate}
\end{theorem}

\subsection{The internal construction}
\margin{sec:rota} \label{sec:rota} We now describe the
construction of a standard Baxter algebra \cite{G-K2},
generalizing Rota~\cite{Ro1}.

For each $n\in \NN_+$, denote $A^{\otimes n}$ for the tensor power
algebra.
Denote the direct limit algebra
\[ \overline{A}= \dirlim A^{\otimes n} = \cup_n A^{\otimes n} \]
where the transition map is given by
\[ A^{\otimes n} \to A^{\otimes (n+1)},\ x \mapsto x\otimes \bfone_A.\]
Note that the multiplication on $A^{\otimes n}$ here is
different from the multiplication on $A^{\otimes n}$ when it
is regarded as the $C$-submodule $\sha_C^{n-1}(A)$ of $\sha_C(A)$.
Let $\frakA(A)$ be the set of sequences with entries in
$\overline{A}$. Thus we have
\[ \frakA(A) = \prod_{n=1}^\infty \overline{A}
    =\left \{ (a_n)_n \mid
    a_n \in \overline{A} \right \}.\]
Define addition, multiplication and scalar multiplication
on $\frakA(A)$ componentwise, making $\frakA(A)$ into a
$\overline{A}$-algebra,
with the all $1$ sequence $(1,1,\ldots)$ as the identity.
Define
\[P_A\tp=P_{A,\lambda}\tp: \frakA(A)\to \frakA(A)\]
by
\[ P_A\tp(a_1,a_2,a_3,\ldots)
=\lambda (0,a_1,a_1+a_2,a_1+a_2+a_3,\ldots).\]
Then $(\frakA(A),P\tp_A)$ is in $\Bax_C$.
For each $a\in A$, define $t^{(a)}=(t^{(a)}_k)_k$ in $\frakA(A)$
by taking
\[ t^{(a)}_k=\otimes_{i=1}^k a_i,\ a_i = \left \{
    \begin{array}{ll} a, & i=k,\\
    1, & i\neq k. \end{array} \right .  \]
Let $\frakS(A)$ be the Baxter subalgebra of $\frakA(A)$ generated
by the sequences
$t^{(a)},\ a\in A$.

\begin{theorem}
\margin{thm:s-r} \label{thm:s-r} \cite{G-K2,Ro1} Assume that the
annihilator of $\lambda\in C$ in the $C$-module $\overline{A}$ is
zero. The morphism in $\Bax_C$
\[ \Phi: \sha_C(A)\to \frakS(A) \]
induced by sending
$a\in A$ to $t^{(a)}$
is an isomorphism.
Therefore,
$(\frakS(A),P_A\tp)$ is a free Baxter algebra on $A$ in the
category $\Bax_C$. \proofend
\end{theorem}

\begin{coro}
\margin{co:rx}
\label{co:rx}
Assume that $\lambda$ is not a zero divisor in $C$.
Let $X$ be a set.
The morphism in $\Bax_C$
\[ \Phi: \sha_C(X)\to \frakS(X)\defeq \frakS(C[X]) \]
induced by sending
$x\in X$ to $t^{(x)}=(t^{(x)}_1,\ldots,t^{(x)}_n,\ldots)$
is an isomorphism.
\end{coro}

There is also an internal construction of free complete Baxter
algebras.

\begin{theorem}
\margin{thm:rcomp}
\label{thm:rcomp}
\cite{G-K2}
Assume that the annihilator of $\lambda\in C$ in
$\overline{A}$ is trivial.
The isomorphism
$\Phi: \sha_C(A)\to \frakS(A)$
extends to an injective homomorphism of Baxter algebras
\[ \widehat{\Phi}: \widehat{\sha}_C(A) \to \frakA(A).\]
\end{theorem}

For $A=C$,
we have
$ \sha_C(C) = \bigoplus_{n\in \NN} C \bfone^{\otimes n} $
and
$ \widehat{\sha}_C(C)= \prod_{n\in\NN} C\bfone^{\otimes n} $
in which the multiplication is given by the equation in
Proposition~\ref{prop:unit}.
Also,
$ \overline{C}= \dirlim C^{\otimes n} \cong C $
and
$ \frakA(C) = \prod_{n=1}^\infty C $
with componentwise addition and multiplication.

\begin{prop}
\margin{prop:src}
\label{prop:src}
Let $C$ be a domain and let $\lambda\in C$ be non-zero.
Then for any
$ b= \sum_{n=0}^\infty b_n \bfone^{\otimes n}\in \sha_C(C)$,
we have
\[ \Phi(b)= \left ( \sum_{i=0}^{n-1}
    \bincc{n-1}{i}\lambda^i b_i \right)_{n\in\NN_+}
    \in \frakS(C).\]
The same formula holds for $\hat{\Phi}$.
\end{prop}
\proof
Since $\Phi$ is $C$-linear, we only to show that, for each
$n\in \NN$,
\margin{eq:src}
\begin{equation}
\Phi(\bfone^{\otimes n})
    = (\bincc{k-1}{n} \lambda^n )_k.
\label{eq:src}
\end{equation}
Note that, by convention, $\bincc{j}{i}=0$ for $j<i$.

We prove equation~(\ref{eq:src}) by induction.
When $n=0$,
$\bfone^{\otimes 0}= \bfone (\defeq \bfone_C)\in C$.
Since $\Phi$ is a $C$-algebra homomorphism,
we have
\[ \Phi(\bfone_C)= (\bfone, \bfone, \ldots)
= (\bincc{k-1}{0} \bfone)_k.     \]
This verifies equation~(\ref{eq:src}) for $n=0$.
Assume that equation~(\ref{eq:src}) is true for $n$.
Then we have
\begin{eqnarray*}
\lefteqn{ \Phi (\bfone^{\otimes (n+1)}) =
    \Phi (P_C(\bfone^{\otimes n}))}\\
&=& P\tp_C (\Phi(\bfone^{\otimes n})) \\
&=& P\tp_C ( (\bincc{k-1}{n} \lambda^n \bfone )_k) \\
&=& \lambda( \sum_{i=1}^{k-1} \bincc{i-1}{n} \lambda^n \bfone)_k \\
&=& ( \bincc{k-1}{n+1} \lambda^{n+1} \bfone)_k.
\end{eqnarray*}
This completes the induction and verifies the first equation
in the proposition.
The second equation follows from the first equation and
Theorem~\ref{thm:rcomp}.
\proofend

\section{Basic properties}
\margin{sec:basic}
\label{sec:basic}
We will first consider subalgebras, quotient algebras
and colimits. Further properties of Baxter algebras will be
studied in later sections.

\subsection{Subalgebras}

\begin{prop}
\margin{prop:sub}
\label{prop:sub}
Let $f: A\rar B$ be an injective $C$-algebra homomorphism,
and let $A$ and $B$ be flat as $C$-modules.
Then the induced Baxter $C$-algebra homomorphisms
$\sha_C(f): \sha_C(A)\rar \sha_C(B)$ and
$\widehat{\sha}_C(f): \widehat{\sha}_C(A)\rar \widehat{\sha}_C(B)$
are injective.
\end{prop}

\proof
By the construction of $\sha_C(A)=\oplus_{n\in \NN_+} A^{\otimes n}$,
$\sha_C(f)$ is defined to be
\[ \bigoplus_{n\in \NN_+} f^{\otimes n}: \bigoplus_{n\in \NN_+} A^{\otimes n}
    \to \bigoplus_{n\in \NN_+} B^{\otimes n}\]
where $f^{\otimes n}: A^{\otimes n}\to B^{\otimes n}$ is
the tensor power of the $C$-module map $f$.
Also by Theorem~\ref{thm:series},
$\widehat{\sha}_C(f)$ can be described as
\[ \prod f^{\otimes n}: \prod_{n\in \NN_+} A^{\otimes n}
    \to \prod_{n\in \NN_+} A^{\otimes n}.\]
Thus we only need to prove that $f^{\otimes n}$ is injective
for all $n\geq 1$.
$f^{\otimes 1}=f$ is injective by assumption.
Assume that $f^{\otimes n}$ is injective.
Since $A$ is flat, $A^{\otimes n}$ is also flat.
So $f:A\rar B$ is injective implies that
\[ \id_{A^{\otimes n}}\otimes f: A^{\otimes (n+1)} = A^{\otimes n}\otimes A
    \rar A^{\otimes n} \otimes B\]
is injective.
By inductive assumption,
$f^{\otimes n}: A^{\otimes n}\rar B^{\otimes n}$
is injective.
Since $B$ is flat,
\[ f^{\otimes n}\otimes \id_B:
A^{\otimes n} \otimes B \rar B^{\otimes n}\otimes B
    =B^{\otimes (n+1)}\]
is injective.
Thus we have that
\[f^{\otimes (n+1)}=(\id_{A^{\otimes n}}\otimes f)
\circ (f^{\otimes n}\otimes \id_B):
A^{\otimes (n+1)} \rar B^{\otimes (n+1)}\]
is injective,
finishing the induction.
\proofend

\subsection{Baxter ideals}
We now study Baxter ideals of $\sha_C(A)$ generated by
ideals of $A$.
Let $I$ be an ideal of $A$. For each $n\in \NN$, let
$I^{(n)}$ be the $C$-submodule of $\sha_C(A)$ generated by
the subset
$ \{\otimes_{i=0}^n x_i \mid\ x_i\in A,\
    x_i\in I {\rm\ for\ some\ } 0\leq i\leq n\}.$

\begin{prop}
\margin{prop:coe}
\label{prop:coe}
Let $I$ be an ideal of $A$.
Let $\tilde{I}$ be the
Baxter ideal of $\sha_C(A)$ generated by $I$
and let $\hat{I}$ be the
Baxter ideal of $\widehat{\sha}_C(A)$ generated by $I$.
Then
\[ \tilde{I} = \bigoplus_{k\in \NN} I^{(n)} \subseteq \sha_C(A)\]
and
\[\hat{I} = \prod_{k\in \NN} I^{(n)} \subseteq
    \widehat{\sha}_C(A)\]
\end{prop}

\proof
Denote
\[ S=\{ \otimes_{i=0}^n x_i \mid x_i\in A,\  0\leq i\leq n,
    {\rm \ and\ } x_i\in I {\rm \ for\ some\ } 0\leq i\leq n,
    n\in \NN \}. \]
Then clearly
\[ \bigoplus_{k\in \NN} I^{(n)} = \sum_{x\in S} C x.\]
So to prove $\tilde{I}\subseteq \oplus_{k\in \NN} I^{(n)}$,
we only need to prove
\[ \tilde{I} \subseteq \sum_{x\in S} C x.\]
Let $J$ denote the sum on the right hand side.
Since clearly $I\subseteq J$,
we only need to prove that $J$ is an Baxter ideal.
Clearly $J$ is a $C$-submodule of $\sha_C(A)$
and is closed under the Baxter operator $P_A$.
For any $x\in S$ and
$y=\otimes _{j=0}^m y_j\in A^{\otimes (m+1)}$,
we have

\[x  y=x_0y_0 \otimes \sum_{(\sigma,T)\in \bs(n,m)}
\lambda^{\mid T\mid}
\sigma ((\otimes_{i=1}^n x_i)\otimes (\otimes_{j=1}^m y_j);T).\]
From the definition of $S$, either $x_0\in I$ or
$x_i\in I$ for some $1\leq i\leq n$.
Thus in each term of the above sum, either
$x_0y_0\in I$ or one of the tensor factors of
$\sigma ((\otimes_{i=1}^n x_i)\otimes (\otimes_{j=1}^m y_j))$
is in $I$. This shows that $x  y\in J$.
Thus $J$ is an Baxter ideal of $\sha_C(A)$.
This proves $\tilde{I}\subseteq \oplus_{k\in \NN} I^{(n)}$.

We next prove by induction on $n$ that each
$I^{(n)}$ is in $\tilde{I}$.
When $n=0$, then $x\in I^{(n)}$ means that $x\in I$.
So the claim is true. Assuming that the claim is true
for $n$ and let
$x=\otimes_{i=0}^{n+1} x_i\in I^{(n+2)}$.
Then one of $x_i,\ 0\leq i\leq n+1$ is in $I$.
If $x_0\in I$, then
$x=x_0  (1\otimes x_1\otimes \ldots\otimes x_{n+1})$
is in $\tilde{I}$
since $\tilde{I}$ is the ideal of $\sha_C(A)$ generated by $I$.
If $x_i\in I$ for some $1\leq i\leq n+1$,
then in
\[x=x_0  (1\otimes x_1\otimes \ldots\otimes x_{n+1})
= x_0  P_A(x_1\otimes \ldots\otimes x_{n+1}), \]
$x_1\otimes \ldots\otimes x_{n+1}\in \tilde{I}$
by induction.
Thus we again have $x\in \tilde{I}$.
Since $\tilde{I}$ is a $C$-submodule, we have
$I^{(n+1)} \subseteq \tilde{I}$.
This completes the induction.
Therefore, $\oplus_{n\in \NN} I^{(n)} \subseteq \tilde{I}.$
This proves the first equation in the proposition.

To prove the second equation, note that by the construction
of the isomorphism
\[\psi_A: \invlim (\sha_C(A)/\Fil^n \sha_C(A))
    \to \prod_{k\in \NN} A^{\otimes (k+1)}\]
in Theorem~\ref{thm:series},
\begin{eqnarray*}
\prod_{k\in \NN} I^{(k)}
&\cong& \invlim (\oplus_{k\in \NN} I^{(k)} +\Fil^n \sha_C(A))/
    \Fil^n \sha_C(A)\\
&=& \invlim (\tilde{I} +\Fil^n \sha_C(A))/\Fil^n \sha_C(A).
\end{eqnarray*}
So we only need to prove that

\[L\defeq \invlim (\tilde{I} +\Fil^n \sha_C(A))/\Fil^n \sha_C(A)\]
is the Baxter ideal $\hat{I}\tp$ of $\invlim (\sha_C(A)/\Fil^n
\sha_C(A))$ generated by $I$. For each $n\in \NN$, $(\tilde{I}
+\Fil^n \sha_C(A))/\Fil^n \sha_C(A)$ is a Baxter ideal of
$\sha_C(A)/\Fil^n \sha_C(A)$. So the inverse limit $\invlim
(\tilde{I} +\Fil^n \sha_C(A))/\Fil^n \sha_C(A)$ is a Baxter ideal
of $\invlim (\sha_C(A)/\Fil^n\sha_C(A))$. Therefore, $\hat{I}\tp
\subseteq L$.

On the other hand, since $\hat{I}\tp$ is a Baxter ideal
of $\invlim (\sha_C(A)/\Fil^n \sha_C(A))$ containing $I$,
its image $\hat{I}\tp_n$
in $\sha_C(A)/\Fil^n \sha_C(A)$ is a Baxter ideal
containing $(I+\Fil^n \sha_C(A))/\Fil^n \sha_C(A)$.
By the same argument as in the proof of the first equation,
we obtain that the Baxter ideal of $\sha_C(A)/\Fil^n\sha_C(A)$
containing $(I+\Fil^n \sha_C(A))/\Fil^n \sha_C(A)$ is
\[( \oplus_{k\leq n} I^{(k)} +\Fil^n \sha_C(A))/\Fil^n\sha_C(A)
= (\tilde{I} +\Fil^n \sha_C(A))/\Fil_n \sha_C(A).\]
Therefore, $\hat{I}\tp_n\supseteq
(\tilde{I} +\Fil^n \sha_C(A))/\Fil_n \sha_C(A)$.
Taking the inverse limit, we obtain
$\hat{I}\tp \supseteq L$, proving the second equation.
\proofend

\subsection{Quotient algebras}
We can now describe how quotients are preserved under
taking free Baxter algebras.

\begin{prop}
\margin{prop:quo}
\label{prop:quo}
Let $I$ be an ideal of $A$.
Let $\tilde{I}$ be the Baxter ideal
of $\sha_C(A)$ generated by $I$
and let $\hat{I}$ be the Baxter ideal
of $\widehat{\sha}_C(A)$ generated by $I$.
Then
\[ \sha_C(A/I) \cong \sha_C(A)/\tilde{I}\]
and
\[ \hat{\sha}_C(A/I) \cong \hat{\sha}_C(A)/\hat{I}\]
as Baxter $C$-algebras.
\end{prop}

\proof
Let $\pi:A\rar A/I$ and
$\tilde{\pi}:\sha_C(A)\rar \sha_C(A)/\tilde{I}$ be the
natural surjections.
The composite map
\[ A\ola{j_A} \sha_C(A)\ola{\tilde{\pi}} \sha_C(A)/\tilde{I}\]
has kernel $I$ by Proposition~\ref{prop:coe}.
Let $j_A\tp:A/I \rar \sha_C(A)/\tilde{I}$
be the induced embedding.
We only need to verify that $\sha_C(A)/\tilde{I}$ with
the Baxter operator $P_A\tp$ induced from $P_A$,
and the embedding $j_A\tp$ satisfies the universal property
for a free Baxter $C$-algebra on $A/I$.

Let $(R,P)$ be an Baxter $C$-algebra and let
$\varphi:A/I \rar R$ be a $C$-algebra homomorphism.
By the universal property of $\sha_C(A)$, the
$C$-algebra homomorphism
\[ \eta\defeq \varphi\circ \pi: A\rar R\]
extends uniquely to an Baxter $C$-algebra homomorphism
\[ \tilde{\eta}: (\sha_C(A),P_A)\rar (R,P).\]
Since $I$ is in the kernel of $\eta$, $\tilde{I}$ is
in the kernel of $\tilde{\eta}$, thus
$\tilde{\eta}$ induces uniquely an Baxter $C$-algebra
homomorphism

\[ \tilde{\eta}\tp: (\sha_C(A),P_A\tp)\rar (R,P). \]
We can summarize these maps in the following diagram

\[\xymatrix{
& A  \ar[dl]_\pi \ar[ddr]^(0.4)\eta\ \ \ar@{^{(}->}[rr]^{j_A}
&&\sha_C(A) \ar[ddl]_(0.4){\tilde{\eta}} \ar[dr]^{\tilde{\pi}} &\\
A/I \ \ \ \
\ar[rrrr]^(0.4){j_A\tp} |!{[u];[d]}\hole |!{[ur];[drr]}\hole
|!{[urrr];[drr]}\hole
\ar[drr]_\varphi
&&&& \ \ \sha_C(A)/\tilde{I} \ar[dll]^{\tilde{\eta}\tp} \\
&& R && } \]

We have
\begin{eqnarray*}
\varphi\circ\pi &=& \eta\ ({\rm\ by\ definition})\\
&=& \tilde{\eta}\circ j_A\
    ({\rm\ by\ freeness\ of\ }\sha_C(A) {\rm \ on\ } A)\\
&=& \tilde{\eta}\tp\circ \tilde{\pi}\circ j_A\ ({\rm\ by\ definition})\\
&=& \tilde{\eta}\tp \circ j_A\tp \circ \pi\ ({\rm\ by\ definition}).
\end{eqnarray*}
Since $\pi$ is surjective, we have
$\varphi=\tilde{\eta}\tp \circ j_A$.
If there is another $\tilde{\eta}\ttp$ such that
$\varphi=\tilde{\eta}\ttp \circ j_A$,
then we have

\begin{eqnarray*}
\tilde{\eta}\tp \circ j_A =\tilde{\eta}\ttp
& \Rightarrow&
\tilde{\eta}\tp \circ j_A \circ \pi =\tilde{\eta}\ttp \circ \pi\\
& \Rightarrow &
\tilde{\eta}\tp \circ \tilde{\pi} \circ j_A =
    \tilde{\eta}\ttp \circ \tilde{\pi} \circ j_A\
    ({\rm\ by\ definition}) \\
& \Rightarrow &
\tilde{\eta}\tp \circ \tilde{\pi} =
    \tilde{\eta}\ttp \circ \tilde{\pi}\
    ({\rm\ by\ freeness\ of\ }\sha_C(A) {\rm \ on\ } A)\\
& \Rightarrow&
\tilde{\eta}\tp = \tilde{\eta}\ttp \
    ({\rm\ by\ surjectivity\ of\ } \tilde{\pi}).
\end{eqnarray*}
This proves the first equation of the proposition.

To prove the second equation, consider the following
commutative diagram

\[
\begin{array}{ccccccccc}
0&\to &\tilde{I}&\to&\sha_C(A)&\ola{\sha_C(\pi)}&\sha_C(A/I)&
    \to&0\\
&& \dar &&\dar &&\dar &&\\
0&\to &\hat{I}&\to&\widehat{\sha}_C(A)&
\ola{\widehat{\sha}_C(\pi)}&\widehat{\sha}_C(A/I)&
    \to&0
\end{array}\]
in which the vertical maps are injective.
${}^{}$ From the first part of the proposition,
the top row is exact.
The desired injectivity of the bottom row is clear,
and the desired surjectivity follows from the definition of
$\widehat{\sha}_C(\pi)$.
Also from the description of $\hat{I}$ in Proposition~\ref{prop:coe},
$\hat{I}\subseteq \ker (\widehat{\sha}_C(\pi))$.
On the other hand,
\begin{eqnarray*}
 (x_n)_n\in \ker (\widehat{\sha}_C(\pi))
 &\Leftrightarrow& (\pi^{\otimes (n+1)}(x_n))_n=0 \\
 &\Leftrightarrow& \pi^{\otimes (n+1)}(x_n)=0, \forall n\geq 0\\
 &\Leftrightarrow& x_n\in \ker (\pi^{\otimes (n+1)}),
    \forall n\geq 0\\
 &\Rightarrow& x_n\in \ker (\sha_C(\pi)), \forall n\geq 0 \\
 &\Rightarrow& x_n\in I^{(n)}, \forall n\geq 0\\
 &\Rightarrow& (x_n)\in \hat{I}.
\end{eqnarray*}
This proves the exactness of the bottom row, hence the
second equation in the proposition.
\proofend

\subsection{Colimits}

\begin{prop}
\margin{prop:lim}
\label{prop:lim}
Let $\Lambda$ be a category whose objects form a set.
Let $F:\Lambda \rar \Alg_C$ be a functor.
Denote $A_\lambda$ for $F(\lambda)$, and
denote $\colim_{\lambda}$ for the
colimit over $\Lambda$.
Then
$\colim_{\lambda} (\sha_C(A_\lambda,P_{A_\lambda}))$
exists and
\[ \colim_{\lambda} (\sha_C(A_\lambda), P_{A_\lambda})
\cong (\sha_C(\colim_{\lambda}A_\lambda),
    P_{\colim_{\lambda}A_\lambda}).\]
In particular, for $C$-algebras $A$ and $B$,
$ \sha_C(A\otimes B) $ is the coproduct of
$\sha_C(A)$ and $\sha_C(B)$.
\end{prop}

\proof
It is well-known that colimits exist in $\Alg_C$.
The proposition then follows from the dual of
~\cite[Theorem 1, p114]{Ma}, stated in page 115.
\proofend

Similar statement for the complete free Baxter algebra
is not true.
For example,
let $\Lambda=\NN_+$ and for each $n\in \Lambda$, let
$A_n=C[x_1,\ldots,x_n]$. With the natural inclusion,
$\{A_n\}$ is a direct system, with
$\colim_n A_n=C[x_1,\ldots,x_n,\ldots]$.
We have
$\colim_n (\widehat{\sha}_C(A_n),P_{A_n})=
\cup_n (\widehat{\sha}_C(A_n),P_{A_n})$
and
$\widehat{\sha}_C(\colim_n A_n)=
    \widehat{\sha}_C(\cup_n A_n)$.
The element $(\otimes_{i=1}^{k+1} x_i)_{k\in \NN}$ is in
$\widehat{\sha}_C(\cup_n A_n)$. But it is not in any
$\widehat{\sha}_C(A_n)$, and hence is not in $\cup_n
(\widehat{\sha}_C(A_n),P_{A_n})$.

\section{Integral domains and reduced algebras}
\margin{sec:domain}
\label{sec:domain}
In this section, we investigate the question of
when a free Baxter $C$-algebra or a free complete Baxter
$C$-algebra is a domain and when  it is
a reduced algebra. We also study the nilpotent elements
when the free Baxter algebra is not reduced.
We will consider the case when $C$ has characteristic zero
in Section~\ref{ss:zero}, and consider the case when $C$
has positive characteristic in Section~\ref{ss:pos}.

\subsection{Case 1: $C$ has characteristic zero}
\margin{ss:zero}
\label{ss:zero}

We begin with the special case when $C$ is a field.
The general case will be reduced to this case.

\subsubsection{$\sha_C(A)$ and $\widehat{\sha}_C(A)$
    when $C$ is a field}

\begin{prop}
\margin{prop:domain1}
\label{prop:domain1}
Let $C$ be a field of characteristic zero.
Assume that $A$ is a $C$-algebra and an integral domain.
\begin{enumerate}
\item
$\sha_C(A)$ is
an integral domain for any $\lambda$.
\item
$\widehat{\sha}_C(A)$ is an integral domain if and only if
$\lambda=0$.
\end{enumerate}
\end{prop}

\proof
1.
Let $\Sigma$ be a basis set of $A$ as a vector space over $C$,
and let $\prec$ be a linear order on $\Sigma$, assuming the
axiom of choice.
Thus
\[A=\bigoplus_{\mu\in \Sigma} C \mu \]
and consequently,
\[ A^{\otimes n} = \bigoplus_{{\mu}\in \Sigma^n}
    C {\mu} \]
where ${\mu}    =(\mu_1,\ldots,\mu_n)\in \Sigma^n$.
Let
\[ \Sigma^\infty =\bigcup_{n\geq 1} \Sigma^n\ {\rm and\ }
\overline{\Sigma}^\infty =\bigcup_{n\geq 0} \Sigma^n\]
with the convention that $\Sigma^0$ is the singleton $\{\phi\}$.
In the following,
we identify a vector $(\mu_1,\ldots,\mu_n)\in \Sigma^n$ with
the corresponding tensor
$\mu_1\otimes \ldots\otimes \mu_n\in A^{\otimes n}$.
Then as a $C$-vector space,
\[\sha_C(A) =\bigoplus_{{\mu}\in \Sigma^\infty}
    C {\mu}. \]
It follows that, as a $A$-module,
\begin{eqnarray*}
 \sha_C(A)
&=&\bigoplus_{{\mu}\in \overline{\Sigma}^\infty}
    A \otimes {\mu}\\
&=&\bigoplus_{{\mu}\in \overline{\Sigma}^\infty}
    A\ (\bfone_A\otimes {\mu})
\end{eqnarray*}
with the convention that $\bfone_A \otimes \phi =\bfone_A$.

We next endow $\overline{\Sigma}^\infty$ with the following
variant of the lexicographic
order induced from the order $\prec$ on $\Sigma$.
We define the empty set $\phi$ to be the smallest element and,
for ${\mu}\in \Sigma^m$ and ${\nu}\in \Sigma^n$,
$m,\ n>0$,
define ${\mu}\prec {\nu}$ if $m<n$, or
$m=n$ and for some $1\leq m_0\leq m$ we have
$\mu_{m_0}\prec \nu_{m_0}$ and  $\mu_i=\nu_i$ for
$m_0+1 \leq i\leq m$.
We also denote this order on $\overline{\Sigma}^\infty$
by $\prec$. It is a linear order.
It is easy to check that, if ${\mu}\prec {\mu}\tp$,
then for any ${\nu}\in \overline{\Sigma}^\infty$,
\margin{eq:order}
\begin{equation}
\label{eq:order} \max \{ \xi \mid \xi \in S({\mu},{\nu})\} \prec
\max \{ \xi \mid \xi \in S({\mu}\tp,{\nu})\}.
\end{equation}
Here
\[ S(\mu,\nu)=\{ \sigma(\mu\otimes \nu) \mid
    \sigma\in S(m,n) \} \]
denotes the set of shuffles
of $\mu$ and $\nu$.
Now let
\[ x=\sum_{{\mu}\in \overline{\Sigma}^\infty}
    a_{{\mu}} (\bfone_A\otimes {\mu}),\
    a_{{\mu}}\in A, \]
and
\[ y=\sum_{{\nu}\in \overline{\Sigma}^\infty}
    b_{{\nu}} (\bfone_A\otimes {\nu}), b_{{\nu}}\in A\]
be two non-zero elements in $\sha_C(A)$.
When $\lambda=0$, only admissible pairs $(\sigma,T)\in \bar{S}(m,n)$
with empty $T$ contribute to the mixable shuffle product
defined in equation~(\ref{eq:shuf}).
So we have
\begin{eqnarray*}
x   y &=&
\sum_{{\mu},{\nu}\in \overline{\Sigma}^\infty}
a_\mu b_\nu \sum_{\xi\in S(\mu,\nu)}
    \bfone_A\otimes \xi\\
&=&
\sum_{\xi \in \overline{\Sigma}^\infty}
c_\xi\ (\bfone_A\otimes \xi).
\end{eqnarray*}
With these notations, we define
\[{\mu}_0
    =\max \{ {\mu} \mid a_{{\mu}}\neq 0\},\]
\[{\nu}_0
    =\max \{ {\nu} \mid b_{{\nu}}\neq 0\}\]
and
\[\xi_0
    =\max \{ \xi \mid c_\xi\neq 0\}.\]
Then from the inequality~(\ref{eq:order}), we have
\[c_{\xi_0} = a_{\mu_0}b_{\nu_0} n_0\]
where $n_0$ is the number of times that $\xi_0$ occurs as a
shuffle of $\mu_0$ and $\nu_0$. If $\lambda\neq 0$, then there are
extra terms in the equation~(\ref{eq:shuf}) of $x  y$ that come
from the mixable shuffles with admissible pairs in which $T$ is
non-empty. But these terms will have shorter lengths and hence are
smaller in the order $\prec$ than the terms from shuffles without
any admissible pairs. So $c_{\xi_0}$ given above is still the
coefficient for the largest term. Note that $n_0$ is a positive
integer by definition. Since $A$ is a domain, we have
$a_{\mu_0}b_{\nu_0}\neq 0$. Since $A$ has characteristic zero, we
further have $a_{\mu_0}b_{\nu_0} n_0\neq 0$. Since
\[ x   y = \sum_{\xi \in \overline{\Sigma}^\infty}
c_\xi\ (\bfone_A\otimes \xi) \]
is the decomposition of $x   y$ according to the
basis $\overline{\Sigma}^\infty$ of the free $A$-module
\[\sha_C(A)= \sum_{\xi\in \overline{\Sigma}^\infty}
    A (\bfone_A\otimes \xi),\]
it follows that
$x   y\neq 0$.

\noindent
2.
Let $\lambda=0$.
The same proof as above, replacing $\max$ by $\min$, shows that
$\widehat{\sha}_C(A)$ is an integral domain.

Let $\lambda\neq 0$.
Consider $x=\bfone_A^{\otimes 2}$,
$y=\sum_{n=0}^\infty (-\lambda)^{-n} \bfone_A^{\otimes (n+1)}$
in $\widehat{\sha}_C(A)$.
By Proposition~\ref{prop:unit}, we have
\begin{eqnarray*}
x   y &=& \sum_{n=0}^\infty
    (-\lambda)^{-n} \bfone_A^{\otimes 2}   \bfone_A^{\otimes (n+1)}\\
&=& \sum_{n=0}^\infty (-\lambda)^{-n}
    \left (\bincc{n+1}{n}\bincc{n}{0} \lambda^0
    \bfone_A^{\otimes (n+2)}
    +\bincc{n}{n}\bincc{n}{1} \lambda \bfone_A^{\otimes (n+1)}
    \right )\\
&=& \sum_{n=0}^\infty (-\lambda)^{-n} (n+1) \bfone_A^{\otimes(n+2)}
+ \sum_{n=0}^\infty (-\lambda)^{-n} n \lambda \bfone_A^{\otimes(n+1)}\\
&=& \sum_{n=0}^\infty (-\lambda)^{-n} (n+1) \bfone_A^{\otimes(n+2)}
- \sum_{n=1}^\infty (-\lambda)^{-n+1} n  \bfone_A^{\otimes(n+1)}\\
&=& 0.
\end{eqnarray*}
So $\widehat{\sha}_C(A)$ is not an integral domain.
\proofend

\subsubsection{$\sha_C(A)$ for a general ring $C$}
Now let $C$ be any ring.
For a $C$-module $N$, denote
\[ N_\tor = \{ x\in N\mid  r x=0
{\rm \ for\ some\ } r\in C,\ r\neq 0\}\] for the $C$-torsion
submodule of $N$. For a domain $D$, denote $\Fr(D)$ for the
quotient field of $D$.

\begin{theorem}
\margin{thm:domain}
\label{thm:domain}
Let $A$ be a $C$-algebra of characteristic zero,
with the $C$-algebra structure given by $\varphi: C\to A$.
Denote $I_0=\ker \varphi$.
The following statements are equivalent.
\begin{enumerate}
\item
$\sha_C(A)$ is a domain.
\item
$A$ is a domain and $(A^{\otimes n})_\tor =I_0,\ $ for all $n\geq 1$.
\item
$A$ is a domain and
the natural map $A^{\otimes n}
    \rar \Fr(C/I_0)\otimes A^{\otimes n}$
is injective for all $n\geq 1$.
\end{enumerate}
\end{theorem}

\proof
Let $\bar{C}=C/I_0$. Then $A$ is also a $\bar{C}$-algebra.
It is well-known that the tensor product $A\otimes_C A$ is
canonically isomorphic to $A\otimes _{\bar{C}}A$ as
$C$-modules and as $\bar{C}$-modules.
It follows that, as a ring,
the $C$-algebra $\sha_C(A)$ is canonically
isomorphic to the $\bar{C}$-algebra $\sha_{\bar{C}}(A)$.
Since being an integral domain is a property of
a ring, $\sha_C(A)$ is a domain if and only if
$\sha_{\bar{C}}(A)$ is one.
Similarly, $\widehat{\sha}_C(A)$ is a domain if and only if
$\widehat{\sha}_{\bar{C}}(A)$ is one.
Thus we only need prove the theorem in the case when
$\varphi:C\to A$ is injective.
So we can assume that $I_0=0$.
We will make this assumption for the rest of the proof.

First note that if $A$ is a domain, then $C$ is also a domain.
In this case we denote $S=C-\{0\}$ and $F=\Fr(C)$.

\noindent
($2 \Leftrightarrow 3$).
This follows from the fact~\cite[Exercise 3.12]{A-M} that,
for each $n\geq 1$,
\[
 (A^{\otimes n})_\tor=
    \ker\{A^{\otimes n} \rar F\otimes_C A^{\otimes n} \}.\]
Therefore the second and the third statement are equivalent.

\noindent
($3 \Rightarrow 1$).
We have the natural isomorphisms
$F\cong S^{-1}C$, $F\otimes A\cong S^{-1} A$
and
\margin{eq:tensor}
\begin{equation}
\label{eq:tensor}
S^{-1}(A^{\otimes n}) \cong
(S^{-1}A)_F^{\otimes n} \defeq
\underbrace{S^{-1}A\otimes_F \ldots \otimes_F S^{-1}A}
    _{n-{\rm factors}}
\cong (S^{-1}A)^{\otimes n}.
\end{equation}
Here the first isomorphism is from~\cite[Proposition 3.3.7]{A-M}
and the last isomorphism follows from the definition of
tensor products and the assumption that $C$ and $A$ are domains.
By the universal property of $\sha_C(A)$
as a free Baxter algebra,
the natural $C$-algebra homomorphism
$f: A\rar S^{-1}A$
gives a $C$-algebra homomorphism
$ \sha_C(f): \sha_C(A)\rar \sha_C(S^{-1}A). $
In fact,
$\sha_C(f)=\oplus_{n=1}^\infty f^{\otimes n}$
where $f^{\otimes n}$ is the tensor power of $f$.
By equation~(\ref{eq:tensor}),
\[f^{\otimes n}: A^{\otimes n} \rar (S^{-1}A)^{\otimes n}
    \cong (S^{-1}A)_F^{\otimes n}
    \cong F\otimes (A^{\otimes n}).\]
Thus we have a $C$-algebra homomorphism
$\tilde{f}: \sha_C(A)\rar \sha_F(S^{-1}A)$
and,  by the third statement of the proposition,
$f^{\otimes n}$ is injective.
Therefore $\sha_C(A)$ is identified with a $C$ subalgebra of
$\sha_F(S^{-1}A)$ via $\tilde{f}$,
and hence is a domain since $\sha_F(S^{-1}A)$ is a domain
by Proposition~\ref{prop:domain1}.

($1 \Rightarrow 2$).
If $\sha_C(A)$ is a domain, then its subring $A$ is a domain.
Since $C$ is a subring of $A$ and hence of $\sha_C(A)$,
we have $\sha_C(A)_\tor=0$.
Since $\sha_C(A)=\oplus_{n\in \NN_+} A^{\otimes n}$,
$(A^{\otimes n})_\tor=0$ for all $n\in \NN_+$.
\proofend

\begin{coro}
\margin{co:domain}
\label{co:domain}
Let $C$ be a domain of characteristic zero.
\begin{enumerate}
\item
If $A$ is a flat $C$-algebra, i.e., $A$ is a $C$-algebra and
is flat as a $C$-module, then $\sha_C(A)$ is a domain.
In particular, for any set $X$, $\sha_C(X)$ is a domain.
\item
If $C$ is a Dedekind domain, then for a $C$-algebra $A$, the free
Baxter algebra $\sha_C(A)$ is a domain if and only if $A$ is
torsion free.
\end{enumerate}
\end{coro}
\proof
As in the proof of Theorem~\ref{thm:domain}, we can assume
that $C$ is a subring of $A$.

\noindent
1.
If $A$ is a flat $C$-module, then $A^{\otimes n},\ n\geq 1$
are flat $C$-modules, so from the injective map
$C\rar F$ of $C$-modules, we obtain the injective map
$A^{\otimes n} =C\otimes A^{\otimes n} \rar F\otimes A^{\otimes n}$.
Hence by Theorem~\ref{thm:domain},
$\sha_C(A)$ is a domain.

\vspace{.3cm}
\noindent
2.
If $C$ is a Dedekind domain, then $A$ is a flat $C$-module
if and only if $A$ is torsion free. Hence the statement.
\proofend

\begin{coro}
\margin{co:prime}
\label{co:prime}
Let $A$ be a $C$ algebra given by the ring homomorphism
$\varphi: C\to A$. Let $I$ be a prime ideal of $A$.
The Baxter ideal $\tilde{I}$ (see Proposition~\ref{prop:coe})
of $\sha_C(A)$
generated by $I$ is a prime ideal if and only if
$A/I$ has characteristic zero and
$((A/I)^{\otimes n})_\tor = \ker \varphi$ for all $n\geq 1$.
\end{coro}

\proof
By Proposition~\ref{prop:quo}, $\tilde{I}$ is a prime ideal
if and only if
$\sha_C(A/I)$ is a domain.
If $A/I$ has characteristic zero and
$((A/I)^{\otimes n})_\tor = \ker \varphi$ for all $n\geq 1$,
then by Theorem~\ref{thm:domain},
$\sha_C(A/I)$ is a domain.
Conversely, if $A/I$ has non-zero characteristic,
then by Theorem~\ref{thm:red} (the proof of which is independent of
Theorem~\ref{thm:domain}), $\sha_C(A/I)$ is not a domain.
If $A/I$ has zero characteristic, but
$((A/I)^{\otimes n})_\tor \neq \ker \varphi$ for some $n\geq 1$,
then by Theorem~\ref{thm:domain},
$\sha_C(A/I)$ is not a domain.
This proves the corollary.
\proofend

\subsubsection{$\widehat{\sha}_C(A)$ for a general ring $C$}
We now consider complete Baxter algebras.

\begin{lemma}
\margin{lem:coe}
\label{lem:coe}
Let $C$ be a UFD and let $x \in \widehat{\sha}_C(X)$ be non-zero.
\begin{enumerate}
\item
Let $\lambda\in C$ be a prime element.
There is
$m\in \NN$ such that $x=\lambda^m x\tp$ and such that
$x\tp \not\in \lambda \widehat{\sha}_C(X)$.
\item
Let $\lambda\in C$ be non-zero.
If $\lambda x=0$, then $x =0$.
\end{enumerate}
\end{lemma}

\proof
1.
By Theorem~\ref{thm:series}, any element $x\in \widehat{\sha}_C(X)$
has a unique expression of the form

\[ x= \sum_{n=0}^\infty x_n,\ x_n\in \sha_C^n(X) = C[X]^{\otimes (n+1)}. \]
So $x\neq 0$ if and only if $x_{n_0}\neq 0 $ for some $n_0\in \NN$.
Let $M(X)$ be the free commutative monoid on $X$.
Define

\[ \overline{X}_n = \{ \otimes_{i=1}^n u_i
    \mid u_i \in M,\ 1\leq i\leq n \} .\]
Then we have

\[ C[X]^{\otimes (n+1)} = \oplus_{u\in \overline{X}_{n+1}} C u \]
and $x_{n_0}$ can be uniquely expressed as
$x_{n_0} = \sum_{u\in \overline{X}_{n_0+1}} c_u u.$
Thus $x_{n_0}\neq 0$ implies that $c_{u_0}\neq 0$ for
some $u_0\in \overline{X}_{n_0+1}$.
Since $C$ is a UFD, there is $m_0\in \NN_+$ such that
$c_{u_0}\not\in \lambda^{m_0} C$.
Then $x_{n_0}\not\in \lambda^{m_0} C[X]^{\otimes (n_0+1)}$
and $x\not\in \lambda^{m_0} \widehat{\sha}_C(X)$.
Therefore the integer
\[ \max \{ k \mid k\in \NN, x\in \lambda^k \widehat{\sha}_C(X) \} \]
exists. This integer can be taken to be the $m$ in the
first statement of the lemma.

\noindent
2.
Assume that $x\in \widehat{\sha}_C(X)$ is non-zero.
Then as in the proof of the first part of the lemma,
there is $n_0\in \NN$ such that
\[x=\sum_{n=0}^\infty x_n,\ x_n\in \sha_C^n(X) = C[X]^{\otimes (n+1)} \]
and $x_{n_0}\neq 0$.
Also, there is $u_0\in \overline{X}_{n_0+1}$ such that

\[ x_{n_0} = \sum_{u\in \overline{X}_{n_0+1}} c_u u,\ c_u\in C\]
and $c_{u_0}\neq 0$.
Since $C$ is a domain and $\lambda \neq 0$,
we have $\lambda c_{u_0}\neq 0$.
Since $C[X]^{\otimes (n_0+1)}$ is a free $C$-module with
the set $\overline{X}_{n_0+1}$ as a basis,
we have
$\lambda x_{n_0+1} \neq 0$.
This in turn proves $\lambda x\neq 0$.
\proofend

\begin{theorem}
\margin{thm:domain1}
\label{thm:domain1}
Let $C$ be a $\QQ$-algebra and a domain
with the property that for every maximal ideal
$M$ of $C$, the localization $C_M$ of $C$ at $M$ is a UFD.
Let $X$ be a set.
For $\lambda\in C$, $\widehat{\sha}_C(X)$ is a domain
if and only if $\lambda$ is not a unit.
\end{theorem}

\noindent {\bf Remarks: } 1. If $C$ is the affine ring of a
nonsingular affine variety on a field of characteristic zero, then
$C$ is locally factorial~\cite[p. 257]{Ei}. Hence
Theorem~\ref{thm:domain1} applies.

\noindent
2. If $C$ is not a $\QQ$-algebra, the statement in the
theorem does not hold. See the example after the proof.

\vspace{0.3cm}

\proof
Assume that $\lambda$ is a unit.
Consider the elements
$x=\bfone_C^{\otimes 2}$,
$y=\sum_{n=0}^\infty (-\lambda)^{-n} \bfone_C^{\otimes (n+1)}$
in $\widehat{\sha}_{C}(X)$.
As in the proof of Proposition~\ref{prop:domain1},
we verify that
$x  y=0$.
So $\widehat{\sha}_C(X)$ has zero divisors and is not a domain.

Now assume that $\lambda\in C$ is not a unit.
We will prove that $\widehat{\sha}_C(X)$ is a domain.
We will carry out the proof in four steps.

{\bf Step 1: }
We first assume that $C$ is a $\QQ$-algebra and a domain,
and assume that $\lambda\in C$ is zero.
We do not assume that, for every maximal
ideal $M$ of $C$, $C_M$ is a UFD.
We clearly have $(C[X]^{\otimes n})_\tor =0$.
Thus from the proof of $2 \Leftrightarrow 3$
and $3 \Rightarrow 1$ in
Theorem~\ref{thm:domain}, the natural map
\[  C[X]^{\otimes n} \to
    (\Fr (C)[X])^{\otimes n}_{\Fr(C)} \]
is injective.
So
\[ \widehat{\sha}_C(X) \to \widehat{\sha}_{\Fr(C)}(X) \]
is injective.
By Proposition~\ref{prop:domain1},
$\widehat{\sha}_{\Fr(C)}(X)$ is a domain.
Therefore,
$\widehat{\sha}_C(X)$ is a domain.

{\bf Step 2: }
Next assume that
$C$ is a UFD and $\lambda\in C$ is a prime element.
Then the ideal $\lambda C$ of $C$ is a prime ideal.
Hence $C/\lambda C$ is a domain.
Then $C[X]/\lambda C[X]\cong (C/\lambda C)[X]$ is
also a domain.
Note that

\begin{eqnarray*}
\lefteqn{ \widehat{\sha}_C(C[X]/\lambda C[X])
 \cong \widehat{\sha}_{C/\lambda C}(C[X]/\lambda C[X])}\\
& \cong & \widehat{\sha}_{C/\lambda C}((C/\lambda C)[X]) \\
&=&\widehat{\sha}_{C/\lambda C}(X)
\end{eqnarray*}
as Baxter $C$-algebras. $C/\lambda C$ is a $\QQ$-algebra and a
domain. So from the first step of the proof, the weight 0 Baxter
$C/\lambda C$-algebra $\widehat{\sha}_{C/\lambda C}(X)$ is a
domain. So $\widehat{\sha}_C(C[X]/\lambda C[X])$ is a domain. On
the other hand, by Proposition~\ref{prop:coe}, $\lambda
\widehat{\sha}_C (X)$ is the Baxter ideal of $\widehat{\sha}_C(X)$
generated by $\lambda C[X]$. So by Proposition~\ref{prop:quo},

\[ \widehat{\sha}_C(C[X]/ \lambda C[X]) \cong
    \widehat{\sha}_C (X)/\lambda \widehat{\sha}_C(X). \]
Thus $\lambda \widehat{\sha}_C (X)$ is a prime ideal.

Suppose $\widehat{\sha}_C(X)$ is not a domain.
Then there are non-zero elements $x,\ y\in \widehat{\sha}_C(X)$
such that
$xy=0$.
From the first statement of Lemma~\ref{lem:coe},
there are $m, n\in \NN$ such that
$ =\lambda^m x\tp,\ y=\lambda^n y\tp$
and $x\tp, y\tp \not\in \lambda\widehat{\sha}_C(X)$.
From the second statement of Lemma~\ref{lem:coe}, we have

\[ xy=0 \Rightarrow \lambda^{m+n}x\tp y\tp =0
    \Rightarrow x\tp y\tp =0.\]
In particular, we have
$x\tp y\tp \in \lambda \widehat{\sha}_C(X)$.
Since $\lambda \widehat{\sha}_C(X)$ is a prime ideal, we must have
$x\tp\in \lambda \widehat{\sha}_C(X)$ or
$y\tp\in \lambda \widehat{\sha}_C(X)$.
This is contradiction.

{\bf Step 3: }
We next assume that
the $\QQ$-algebra $C$ is a domain and a UFD,
and assume that $\lambda\in C$ is not a unit.
Then there is a prime element $\lambda_1\in C$ such that
$\lambda=\lambda_1 \lambda_2$ for some $\lambda_2\in C$.
From the second step of the proof, the weight $\lambda_1$ complete
shuffle algebra
$(\widehat{\sha}_{C,\lambda_1}(X), P_{X,\lambda_1})$ is a domain.

Define another operator $Q$ on the $C$-algebra
$\widehat{\sha}_{C,\lambda_1}(X)$ by

\[ Q(x)=\lambda_2 P_{X,\lambda_1}(x),\
    x\in \widehat{\sha}_{C,\lambda_1}(X). \]
Then
\begin{eqnarray*}
\lefteqn{ Q(x)Q(y)
= \lambda_2^2 P_{X,\lambda_1}(x)P_{X,\lambda_1}(y)}\\
&=& \lambda_2^2 (P_{X,\lambda_1}(xP_{X,\lambda_1}(y))
    +P_{X,\lambda_1}(yP_{X,\lambda_1}(x))
    +\lambda_1 P_{X,\lambda_1}(xy))\\
&=& \lambda_2 P_{X,\lambda_1}(x \lambda_2 P_{X,\lambda_1}(y))
+ \lambda_2 P_{X,\lambda_1}(y \lambda_2 P_{X,\lambda_1}(x))
+ \lambda\lambda_2 P_{X,\lambda_1}(xy)\\
&=& Q(xQ(y))+Q(yQ(x))+\lambda Q(xy)
\end{eqnarray*}
So $(\widehat{\sha}_{C,\lambda_1}(X),Q)$ is a Baxter algebra of
weight $\lambda$.
Since $(\widehat{\sha}_{C,\lambda}(X),P_{X,\lambda})$ is a free
complete Baxter algebra of weight $\lambda$,
there is a unique homomorphism of weight $\lambda$ complete Baxter
algebras
\[ f: (\widehat{\sha}_{C,\lambda}(X),P_{X,\lambda}) \to
    (\widehat{\sha}_{C,\lambda_1}(X),Q)\]
that extends the identity map on $X$.

\begin{lemma}
\margin{lem:inj}
\label{lem:inj}
The map $f$ is injective.
\end{lemma}
\proof
We first prove that, for $n\in\NN$ and $x\in \sha_{C,\lambda}^n(X)$,
\margin{eq:inj}
\begin{equation}
f(x) = \lambda_2^n x\in \sha_{C,\lambda_1}^n (X).
\label{eq:inj}
\end{equation}
Recall that, as a $C$-module,
$\sha_C^n(X)= C[X]^{\otimes (n+1)}$.
The identify map on $X$ induces the identity map on $C[X]$.
This proves equation~(\ref{eq:inj}) for $n=0$.
Assuming that, for $x\in \sha_{C,\lambda}^n(X)$,
$f(x)=\lambda_2^n x \in \sha_{C,\lambda_1}^n(X)$,
and letting
$x_0\otimes \ldots \otimes x_{n+1}\in \sha_{C,\lambda}^{n+1}(X)$,
we have
\begin{eqnarray*}
\lefteqn{f(x)
    = f(x_0 P_{X,\lambda}(x_1\otimes\ldots\otimes x_{n+1}))}\\
&=& f(x_0) f(P_{X,\lambda}(x_1\otimes \ldots \otimes x_{n+1}))\\
&=& x_0 Q (f(x_1\otimes\ldots\otimes x_{n+1})) \\
&=& x_0 Q(\lambda_2^n x_1\otimes \ldots\otimes x_{n+1})\\
&=& x_0  \lambda_2 P_{X,\lambda_1}(\lambda_2^n x_1\otimes \ldots
    \otimes x_{n+1})\\
&=& \lambda_2^{n+1} x_0 \otimes x_1\otimes \ldots \otimes x_{n+1}.
\end{eqnarray*}
This implies that, for any $y\in \sha_{C,\lambda}^{n+1} (X)$,
we have $f(x)=\lambda^{n+1} x\in \sha_{C,\lambda_1}^{n+1} (X)$.
This completes the induction.

By equation~(\ref{eq:inj}), the restriction of $f$ to
$\sha_{C,\lambda}^n(X)$ gives a map
\[f_n: \sha_{C,\lambda}^n(X) \to \sha_{C,\lambda_1}^n (X)\]
and, from Theorem~\ref{thm:series}, we have
\[ f= \prod_{n\in \NN} f_n. \]
Thus in order to prove the lemma, we only need to prove that
$f_n$ is injective for each $n\in\NN$.
This follows from equation~(\ref{eq:inj}) and Lemma~\ref{lem:coe}.
\proofend

We continue the proof of Theorem~\ref{thm:domain}.
Since $\lambda_1$ is a prime element of $C$, by the first
part of the proof,
$\widehat{\sha}_{C,\lambda_1}(X)$ is a domain.
By Lemma~\ref{lem:inj},
$\widehat{\sha}_{C,\lambda}(X)$ is isomorphic to a subalgebra
of $\widehat{\sha}_{C,\lambda_1}(X)$, hence is also a
domain.

{\bf Step 4:} We now consider the general case when
the $\QQ$-algebra $C$ is a domain
whose localization at each maximal ideal is a UFD.
Let $\lambda$ be a non-unit.
Then there is a maximal ideal $M$ of $C$ containing $\lambda$.
By assumption, $C_M$ is a UFD.
If we regard $C$ as a subring of $C_M$ by the natural
embedding $C\to C_M$ given by the localization map,
then $\lambda\in C$ remains a non-unit in $C_M$.
From the proof in Step 3,
$\widehat{\sha}_{C_M,\lambda}(X)$ is domain.
On the other hand, the natural embeddings
\[ C[X] \hookrightarrow C_M[X] \hookrightarrow \Fr(C)[X] \]
induce the natural morphisms of $C$-modules.
\[ C[X]^{\otimes n} \rightarrow C_M[X]^{\otimes n}
\rightarrow \Fr(C)[X]^{\otimes n}.  \]
Since clearly $(A^{\otimes n})_\tor = 0$ for $A=C[X]$,
by the equivalence $(2)\Leftrightarrow (3)$ in
Theorem~\ref{thm:domain},
the composite of the above two maps is injective.
Therefore the map
$C[X]^{\otimes n} \to C_M[X]^{\otimes n}$
is injective.
Thus
$\widehat{\sha}_{C,\lambda}(X)\to \widehat{\sha}_{C_M,\lambda}(X)$
is injective.
This proves that $\widehat{\sha}_{C,\lambda}(X)$ is a domain.
\proofend

When $C$ is not a $\QQ$-algebra, the situation is much more
complicated and will be the subject of a further study.
Here we just give an example to show that Theorem~\ref{thm:domain1}
does not hold without the assumption that $C$ is a $\QQ$-algebra.

\vspace{.3cm}
\noindent
{\bf Example: }
Let $C=\ZZ$ and $\lambda=2$.
Then $C$ is not a $\QQ$-algebra. But all other conditions
in Theorem~\ref{thm:domain1} are satisfied.
Consider the two elements
\[ x=\sum_{n=1}^\infty (-1)^{n+1} \bfone^{\otimes (n+1)}  \]
and
\[ y= 2+\sum_{n=1}^\infty (-1)^n \bfone^{\otimes (n+1)}. \]
in $\widehat{\sha}_C(X)=\widehat{\sha}_C(C[X])$.
By Proposition~\ref{prop:src}, the $n$-th component of $\Phi(x)$ is
\begin{eqnarray*}
\lefteqn{ \sum_{i=1}^n \bincc{n}{i} (-1)^{i+1} 2^i =
    -\sum_{i=1}^n \bincc{n}{i} (-2)^i}\\
&=& - ((1-2)^n-1) \\
&=& \left \{ \begin{array}{ll} 0, & n {\rm\ even }, \\
    2, & n {\rm\ odd\ } \end{array} \right .
\end{eqnarray*}
and the $n$-th component of $\Phi(y)$ is
\begin{eqnarray*}
\lefteqn{ 2+ \sum_{i=1}^n \bincc{n}{i} (-1)^i 2^i
    = 2+\sum_{i=1}^n \bincc{n}{i}(-2)^i }\\
&=& 2+ ((1-2)^n -1) \\
&=& \left \{ \begin{array}{ll} 2, & n {\rm\ even }, \\
    0, & n {\rm\ odd\ } \end{array} \right .
\end{eqnarray*}
Thus $\Phi(xy)=\Phi(x)\Phi(y)=0$.
Therefore, $xy=0$ by Theorem~\ref{thm:s-r}.
This shows that $\widehat{\sha}_C(X)$ is not a domain.

\subsection{Case 2: $C$ has positive characteristic}
\margin{ss:pos}
\label{ss:pos}
We now consider the case when the characteristic of
$C$ is positive.

\begin{theorem}
\margin{thm:red}
\label{thm:red}
Let $C$ be a ring of positive characteristic and let
$A\supseteq C$ be a $C$-algebra.
\begin{enumerate}
\item
$\sha_C(A)$ (hence $\widehat{\sha}_C(A)$)
is not an integral domain.
\item
If $\lambda=0$, then $\sha_C(A)$ (hence $\widehat{\sha}_C(A)$)
is not reduced.
The nil radical of $\sha_C(A)$ (resp. of $\widehat{\sha}_C(A)$)
is given by
\[ N(\sha_C(A))= N(A) \oplus
    (\bigoplus_{n\in \NN_+}A^{\otimes (n+1)})\]
\[({\rm resp.\ }\ \  N(\widehat{\sha}_C(A))\subseteq N(A) \times
    \prod_{n\in \NN_+}A^{\otimes (n+1)}).\]
\item
If $\lambda\neq 0$, and if, for every $k\geq 1$,
$\lambda$ has trivial annihilator
in the $C$-module $A^{\otimes k}$ and the tensor power algebra
$A^{\otimes k}$ is reduced,
then $\sha _C(A)$ and $\widehat{\sha}_C(A)$ are reduced.
\end{enumerate}
\end{theorem}

\proof
1.
Since $A$ is a subalgebra of $\sha_C(A)$,
it is clear that if $A$ is not a domain, then
$\sha_C(A)$ is not a domain.
So we will assume that $A$ is a domain,
hence the characteristic of $A$ is a prime number $p$.

If $\lambda=0$, then from Proposition~\ref{prop:unit},
$(\bfone\otimes \bfone) ^p=(p!) \bfone^{\otimes p} =0$ in
$\sha_C(C)$. Then $(\bfone_A\otimes \bfone_A)^p =(p!)
\bfone_A^{\otimes (p+1)}=0$ in $\sha_C(A)$. So $\bfone_A\otimes
\bfone_A$ is a zero divisor and $\sha_C(A)$ is not a domain.

We now assume $\lambda \neq 0$.
We first let $A=C$.
Then by Proposition~\ref{prop:src}
the isomorphism of Baxter algebras
\[  \Phi: \sha_C(C) \to \frakS(C)\subseteq \prod_{k\in \NN_+} C \]
sends $\bfone$ to $(\bfone, \bfone, \ldots )$
and sends $\bfone^{\otimes 2}=P_C(\bfone)$ to
\[ P\tp_C(\bfone, \bfone, \ldots )
    = \lambda(0, \bfone, 2\bfone, \ldots).\]
Thus for each $0\leq i\leq p-1$, we have

\[ \Phi(i\lambda \bfone + \bfone^{\otimes 2})
    =\lambda (i\bfone, (1+i)\bfone, (2+i)\bfone, \ldots )\]
So the $n$-th component of $\Phi(i\lambda\bfone +\bfone^{\otimes 2})$
is zero for $n\equiv i \pmod{p}$.
Since the product in $\frakS(C)$ is componentwise,
it follows that the $n$-th component of
$ \prod_{i=0}^{p-1} \Phi(i\lambda \bfone +\bfone^{\otimes 2})$
is zero for all $n$.
Therefore
$ \prod_{i=0}^{p-1} \Phi(i\lambda \bfone +\bfone^{\otimes 2})=0.$
Since $\Phi$ is an algebra isomorphism, we have

\[ \prod_{i=0}^{p-1} (i\lambda \bfone +\bfone^{\otimes 2})=0.\]
Clearly none of $i\lambda \bfone +\bfone^{\otimes 2},\ 0\leq i\leq p-1$,
is zero. So these elements are zero divisors
and $\sha_C(C)$ is not a domain.
For a general $C$-algebra $A$ of characteristic $p$,
the ring homomorphism $\varphi: C\to A$ that defines the $C$-algebra
structure on $A$ induces a Baxter algebra homomorphism
$\sha_C(\varphi): \sha_C(C) \to \sha_C(A)$ sending
$\bfone^{\otimes k}$ to $\bfone_A^{\otimes k}$.
So $\sha_C(\varphi)$ sends
\[ \prod_{i=0}^{p-1} (i\lambda \bfone +\bfone^{\otimes 2})=0\]
to
\[ \prod_{i=0}^{p-1} (i\lambda \bfone_A +\bfone_A^{\otimes 2})=0.\]
This shows that $\sha_C(A)$ is not a domain.

\vspace{0.3cm}

\noindent 2. Let $q>0$ be the characteristic of $C$. If
$\lambda=0$, then as in the proof of the first statement of the
theorem, $(\bfone_A\otimes \bfone_A)^q =q! \bfone_A^{\otimes
(q+1)}=0$ in $\sha_C(A)$. So $\bfone_A\otimes \bfone_A$ is a
nilpotent and $\sha_C(A)$ is not reduced.

Before describing the nil radical, we need some
preparation.
Let $(R,P)$ be a Baxter algebra.
For any $x\in R$, denote $P_x(y)=P(xy)$, $y\in R$.
For any $n\in \NN$, denote
$P_x^n=\underbrace{P_x\circ \ldots\circ P_x}_{n-{\rm times}}$
with the convention that $P_x^0=\id_R$.

\begin{lemma}
\label{lem:power}
\margin{lem:power}
Let $(R,P)$ be a Baxter $C$-algebra of weight zero.
\begin{enumerate}
\item
For $n\in \NN$, $P_x^n(\bfone_R)P_x(\bfone_R)
    =(n+1)P^{n+1}_x(\bfone_R)$.
\item
For $n\in \NN$,
$P(x)^n=n!P_x^n(\bfone_R)$.
\end{enumerate}
\end{lemma}
\proof
We prove both statements by induction on $n$.
The first statement is clearly true for $n=0$.
Assume that it is true for $n$.
Then
\begin{eqnarray*}
P_x^{n+1}(\bfone_R)P_x(\bfone_R) &=&
P(xP^n_x(\bfone_R))P(x) \\
&=& P(xP(xP_x^n(\bfone_R)))+ P(xP_x^n(\bfone_R)P(x))\\
&=& P_x^{n+2}(\bfone_R) + P(x(n+1)P_x^{n+1}(\bfone_R)\\
&=& (n+2)P_x^{n+2}(\bfone_R),
\end{eqnarray*}
completing the induction.

The second statement is again clear for $n=0$.
Assume that the statement if true for $n$.
Then by the first statement of the lemma,
\[ P(x)^{n+1}=P(x)^n P(x)=n! P_x^n(\bfone_R) P_x(\bfone_R)
 =(n+1)! P_x^{n+1}(\bfone_R), \]
completing the induction.
\proofend

Now back to the proof of Theorem~\ref{thm:red}.
We first prove
\[ N(\sha_C(A))\supseteq  N(A) \oplus
    (\bigoplus_{n\in \NN_+}A^{\otimes (n+1)}).\]
Clearly $N(A)\subseteq N(\sha_C(A))$.
Let $n\in \NN_+$ and let
$x=x_0\otimes \ldots \otimes x_n\in A^{\otimes (n+1)}$.
Denote $x=x_0\otimes x^+$, with
$x^+= x_1\otimes \ldots \otimes x_n$.
We have $x=x_0\shpr P(x^+)$.
Then from Lemma~\ref{lem:power},
\[ x^q=x_0^q \shpr P(x^+)^q =
    x_0^q \shpr (q! P^q_{x^+}(\bfone_A))
=0. \] Since any element of $\oplus_{k\in \NN_+}A^{\otimes (k+1)}$
is a finite sum of elements of the form $x_0\otimes \ldots \otimes
x_n$, $n\in \NN^+$, and a finite sum of nilpotent elements is
still nilpotent, we have $\oplus_{k\in \NN^+} A^{\otimes (k+1)}
\subseteq N(\sha_C(A)).$ Therefore,
\[ N(\sha_C(A))\supseteq  N(A) \oplus
    (\bigoplus_{n\in \NN_+}A^{\otimes (n+1)}).\]

Now let $x\in \sha_C(A)$ be nilpotent. $x$ can be uniquely
expressed as $\sum_{k\in \NN} x_k$ with $x_k\in A^{\otimes
(k+1)}$. Since $x$ is nilpotent, $x^n=0$ for some $n\in \NN^+$. By
the definition of the mixed shuffle product $\shpr$ in
$\sha_C(A)$, defined in Eq(\ref{eq:shuf}), we can uniquely express
$x^n$ as $\sum_{k\in \NN} y_k$ with $y_k\in A^{\otimes (k+1)}$ and
$y_0=x_0^n$. By the uniqueness of $y_0$, we have $y_0=0$. So $x_0$
is nilpotent. This shows that $x$ is in $N(A)\oplus
(\bigoplus_{n\in \NN^+} A^{\otimes (n+1)})$.

By the same argument as in the last paragraph, we also obtain
\[ N(\widehat{\sha}_C(A)) \subseteq N(A)\times
(\prod_{n\in \NN^+} A^{\otimes (n+1)}). \] This completes the
proof of 2.

 \vspace{.3cm}

\noindent
3.
We first make a general remark on the mixable shuffle product.
For any two tensors
$a_0\otimes\ldots \otimes a_{m}$ and
$b_0\otimes\ldots \otimes b_{m}$ of $A^{\otimes (m+1)}$,
by the definition of the shuffle product in
$\sha_C(A)$, we can write
\[  (a_0\otimes \ldots \otimes a_m)\shpr
(b_0\otimes\ldots \otimes b_{m})
 =\sum_{i\geq m} x_i, \ x_i\in A^{\otimes (i+1)}, \]
and
\[ x_m= \lambda^m a_0 b_0 \otimes \ldots \otimes a_m b_m. \]
In other words,
\[ x_m=\lambda^m (a_0\otimes \ldots \otimes a_m) \cdot
(b_0\otimes\ldots \otimes b_{m}),\] where we use $\cdot$ to denote
the product in the tensor product algebra $A^{\otimes (m+1)}$. By
the biadditivity of the multiplication in $\sha_C(A)$ and the
multiplication $\cdot$ in $A^{\otimes (m+1)}$, we see that for any
non-zero elements $a,\ b\in A^{\otimes (m+1)}$, the term of
$a\shpr b\in \sha_C(A)$ with degree $m$ equals $\lambda^m (a\cdot
b)$.

Now let $x$ be a non-zero element of $\sha_C(A)$.
Express $x$ as
\[ x=\sum_{i=m}^{n} x_i,\ x_i\in A^{\otimes (i+1)}, 0\leq m\leq i\leq n \]
with $x_m\neq 0$ and $x_n\neq 0$.
It follows from the above remark and the induction that,
for any $k\geq 1$,
if we express
\[ x ^k = \sum_{i\geq m} y_i,\ y_i\in
A^{\otimes (i+1)}, i\geq m, \] then $y_m=\lambda^{m(k-1)}
x_m^{\bullet k}$. Here $x_m^{\bullet k}$ stands for the $k$-th
power of $x_m$ in the tensor product algebra $A^{\otimes (m+1)}$.
If $\lambda$ does not annihilate any non-zero elements in
$A^{\otimes (m+1)}$, then $\lambda^m x_m\neq 0$. If $A^{\otimes
(m+1)}$ is reduced, then we further have $(\lambda^m x_m)^{\bullet
k}\neq 0$. Since $(\lambda^m x_m)^{\bullet k} = \lambda^m y_m$, we
have $y_m\neq 0$. Then $x^k$ is not zero, proving that $\sha_C(A)$
has no non-zero nilpotent elements, hence is reduced.

The same argument can be applied to $\widehat{\sha}_C(A)$,
proving that $\widehat{\sha}_C(A)$ is reduced.
\proofend

\noindent {\bf Acknowledgements. }
The author thanks William Keigher for helpful discussions.

\addcontentsline{toc}{section}{\numberline {}References}

\end{document}